\documentclass[graybox]{svmult}

\usepackage{mathptmx}       
\usepackage{helvet}         
\usepackage{courier}        
\usepackage{type1cm}        
\usepackage{makeidx}         
\usepackage{graphicx}        
\usepackage{multicol}        
\usepackage[bottom]{footmisc}

\makeindex             

\begin{document}

\title*{Skewness-kurtosis adjusted confidence estimators
and significance tests}
\author{Wolf-Dieter Richter}
\institute{Wolf-Dieter Richter \at Institute of Mathematics, University of Rostock,
Ulmenstra{\ss}e 69, Haus 3, 18057 Rostock, Germany \email{name@email.address}
}
%
%
\maketitle

\abstract{First and second kind  modifications of usual confidence intervals for estimating the expectation and of usual local alternative parameter choices  are introduced in a way such that the asymptotic behavior of the true non-covering probabilities and the covering probabilities under the modified local non-true parameter assumption  can be asymptotically exactly controlled.
The orders of convergence to zero of both types of  probabilities are assumed  to be suitably bounded below according to an Osipov-type condition and the sample distribution is assumed to satisfy a corresponding tail condition due to Linnik.
Analogous considerations are presented for the power function when testing a hypothesis concerning the expectation both under the assumption of a true hypothesis as well as under a modified local alternative.  Applications are given for exponential families.}

\section{Introduction}
\label{sec:1}
Asymptotic normality of the distribution of the suitably centered and normalized arithmetic mean of  i.i.d. random variables is one of the best studied and most often exploited facts in asymptotic statistics.
It is supplemented in  local asymptotic normality theory by  limit theorems for the corresponding distributions under the assumption that the mean is shifted of order $n^{-1/2}$. There are many successful simulations and real applications of both types of central limit theorems, and one may ask for  a more detailed explanation of those success. The present note is aimed to present such additional theoretical  explanation under certain circumstances. Moreover, the note is aimed to further stimulate analogous consideration in more general situations and to stimulate checking the new results by simulation. Moreover, based upon the results presented here, it might become attractive to search for  additional explanation to various known simulation results in the area of asymptotic normality which is, however, behind the scope of the present note.

Based upon a large deviation result in \cite{Li}, skewness-kurtosis modifications of usual confidence intervals for estimating the expectation and of usual local alternative parameter choices  are introduced here in a way such that the asymptotic behavior of the true non-covering probabilities and the covering probabilities under the modified local non-true parameter assumption  can be exactly controlled.
The orders of convergence to zero of both types of  probabilities are suitably bounded below by assuming an Osipov-type condition, see \cite{Os}, and the sample distribution is assumed to satisfy a corresponding Linnik condition.

Analogous considerations are presented for the power function when testing a hypothesis concerning the expectation both under the assumption of a true hypothesis as well as under a local alternative.
Finally, applications are given for exponential families.

A  concrete situation where the results of this paper apply is the case sensitive preparing of the machine settings of a machine tool. In this case, second and higher order moments of the manipulated variable do not change from one adjustment to another one and may be considered to be known over time.

It might be another aspect of stimulating further research if one asks for the derivation of limit theorems  in the future being close to those in \cite{Li} but where higher order moments are estimated.

Let $X_1,...,X_n$ be i.i.d. random variables with the common distribution law from a shift family of distributions, $   P_\mu=P(\cdot-\mu) $, where  the expectation equals
 $\mu,$
  $ \mu \in R$, and the variance is  $\sigma^2$.
It is well known that
$T_n=\sqrt{n}(\bar X_n -\mu)/\sigma$
is asymptotically standard normally distributed,
$T_n\sim AN(0,1)$. Hence,  $P_\mu(T_n>z_{1-\alpha})\rightarrow \alpha$,
and
under the
local non-true parameter assumption, $\mu_{1,n}=\mu+\frac{\sigma}
{\sqrt{n}}(z_{1-\alpha}-z_\beta)$,
i.e. if one assumes that a sample is drawn with a shift of location (or with an error in the variable), then
$
P_{\mu_{1,n}}(T_n \leq z_{1-\alpha})= P_{\mu_{1,n}}(\sqrt{n}\frac{\bar X_n-\mu_{1,n}}{\sigma} \leq z_{\beta})\rightarrow \beta$ as $n\rightarrow \infty,$
where     $z_q$ denotes the quantile of order $q$ of the standard Gaussian distribution.

Let $ACI^u= [\bar X_n - \frac{ \sigma}{\sqrt{n}}z_{1-\alpha}, \infty)$
denote the upper asymptotic  confidence interval for $\mu$ where the true non-covering probabilities satisfy the asymptotic relation
\[
 P_\mu(ACI^u {\; does\; not\; cover\; }  \mu)\rightarrow \alpha,\; n\rightarrow \infty.\]
Because $
P_{\mu_{1,n}}(\bar X_n-\frac{\sigma}{\sqrt{n}}z_{1-\alpha}<\mu)
=
P_{\mu_{1,n}}(\sqrt{n}\frac{\bar X_n-\mu_{1,n}}{\sigma}\leq z_\beta)
$, the covering probabilities under $n^{-1/2}$-locally chosen  non-true parameters satisfy
\[
P_{\mu_{1,n}}(ACI^u {\; covers\; }  \mu) \rightarrow \beta, \; n\rightarrow \infty.\]
The aim of this note is to prove refinements of the latter two asymptotic relations where $\alpha=\alpha(n)\rightarrow 0 $ and $\beta=\beta(n)\rightarrow 0$ as $n\rightarrow \infty$, and to prove similar results for two-sided confidence intervals and for the power function when testing corresponding hypotheses.

\section{Expectation estimation}
\subsection{First and second kind adjusted one-sided confidence intervals}
\label{sec:2}
According to \cite{Li}, it is said that a random variable $X$ satisfies the Linnik condition of order $\gamma, \gamma>0, $ if
\begin{equation}
{E}_\mu
\exp\{|X-\mu|^{\frac{4\gamma}{2\gamma+1}}\}
<\infty.
\label{Lgamma}\end{equation}
Let us define the first kind (or first order) adjusted  asymptotic Gaussian quantile  by
\[
z_{1-\alpha(n)}(1)=z_{1-\alpha(n)}
+\frac{g_1}{6\sqrt{n}}
z^2_{1-\alpha(n)}\]
where  $g_1=\it E (X-\it E(X))^3/\sigma^{3/2}$ is the skewness of $X$. Moreover, let the first kind (order)  adjusted upper asymptotic  confidence interval for $\mu$ be defined by
\[
ACI^u(1)=[\bar X_n -\frac{\sigma}{\sqrt{n}}
z_{1-\alpha(n)}(1), \infty)\]
and denote a first kind modified
non-true local parameter choice
by
\[
\mu_{1,n}(1)=\mu_{1,n}+\frac{\sigma g_1}
{6n} (z^2_{1-\alpha(n)}-z^2_{\beta(n)}).\]
 Let us say that the probabilities $\alpha(n) $ and $\beta(n)$ satisfy an Osipov-type  condition of order $\gamma$ if
\begin{equation}
  n^{\gamma}\exp\{\frac{n^{2\gamma}}{2}\}
\cdot \min\{\alpha(n),\beta(n)\}\rightarrow \infty ,\; n\rightarrow \infty.
\label{Os1}
\end{equation}
This condition means that neither $\alpha(n)$ nor $\beta(n) $ tend to zero as fast as or even faster than $n^{-\gamma}\exp\{-n^{2\gamma}/2\},$
i.e. $\min\{\alpha(n), \beta(n)\}\gg n^{-\gamma}\exp\{-n^{2\gamma}/2\}
,$
and that $ \max\{z_{1-\alpha(n)},z_{1-\beta(n)}\}=
o(n^\gamma), n \rightarrow \infty.$

If two functions $f,g$ satisfy the  relation $\lim\limits_{n\rightarrow\infty }f(n)/g(n)=1$
then this asymptotic equivalence will be expressed   as
$f(n)\sim g(n), n\rightarrow \infty.$
\begin{theorem}
 If $\alpha(n)\downarrow 0$, $\beta(n)\downarrow 0$ as $n\rightarrow \infty$ and conditions ($\ref{Lgamma} $) and ($\ref{Os1}$)   are satisfied for $\gamma\in (\frac{1}{6},\frac{1}{4}]$ then
 \[
 P_\mu (ACI^u(1) {\; does\; not\; cover\; } \mu)\sim \alpha(n), \, n\rightarrow \infty
 \]
 and
 \[
 P_{\mu_{1,n}(1)} (ACI^u(1) {\;  covers\; } \mu)\sim \beta(n), \, n\rightarrow \infty.
 \]
\end{theorem}
Let us define the second kind adjusted  asymptotic Gaussian quantile
\[
z_{1-\alpha(n)}(2)=z_{1-\alpha(n)}(1)
+\frac{3g_2-4g_1^2}{72n}
z^3_{1-\alpha(n)}\]
where  $g_2=\it E (X-\it E(X))^4/\sigma^{4}-3$ is the kurtosis of $X$,
the second kind  adjusted upper asymptotic  confidence interval for $\mu$
\[
ACI^u(2)=[\bar X_n -\frac{\sigma}{\sqrt{n}}
z_{1-\alpha(n)}(2), \infty),\]
and a second kind modified
non-true  local parameter choice
\[
\mu_{1,n}(2)
=\mu_{1,n}(1)+
\frac{\sigma(3g_2-4g_1^2)}
{72n^{3/2}}(z^3_{1-\alpha(n)}-z_{\beta(n)}^3)
.\]
\begin{theorem}
 If $\alpha(n)\downarrow 0$, $\beta(n)\downarrow 0$ as $n\rightarrow \infty$ and conditions ($\ref{Lgamma} $) and ($\ref{Os1}$)  are satisfied for $\gamma\in (\frac{1}{4},\frac{3}{10}]$ then
 \[
 P_\mu (ACI^u(2) {\; does\; not\; cover\; } \mu)\sim \alpha(n), \, n\rightarrow \infty
 \]
 and
 \[
 P_{\mu_{1,n}(2)} (ACI^u(2) {\;  covers\; } \mu)\sim \beta(n), \, n\rightarrow \infty.
 \]
\end{theorem}

\begin{remark}
Under the same assumptions, analogous results are true for lower asymptotic   confidence intervals, i.e. for  $ACI^l(s)=(-\infty, \bar X_n+\frac{\sigma}{\sqrt{n}}z^{-}_{1-\alpha}(s)), s=1,2:$
\[
P_\mu(ACI^l(s)\; does\; not\; cover\; \mu)
\sim \alpha(n)
\]
and
\[
P_{\mu^{-}_{1,n}(s)}(ACI^l(s)\; covers\; \mu)
\sim \beta(n),\, n\rightarrow \infty .\]
Here, $z^{-}_{1-\alpha}(s)$ means the quantity $z_{1-\alpha}(s)$ where $g_1$ is replaced by $-g_1, s=1,2,$
and \[
\mu^{-}_{1,n}(s)=\mu-\frac{\sigma}{\sqrt{n}}
(z_{1-\alpha}-z_\beta)+\frac{\sigma g_1}{6n}
(z^2_{1-\alpha}-z^2_\beta)-
\frac{\sigma(3g_2-4g_1^2)}
{72n^{3/2}}(z^3_{1-\alpha}-
z^3_\beta)I_{\{2\}}(s).\]

\end{remark}

\begin{remark}
In many situations where limit theorems are considered as they were in Section 1, the additional assumptions $(\ref{Lgamma})$ and
$(\ref{Os1})$ may, possibly unnoticed, be fulfilled. In such situations, Theorems 1 and 2, together with the following theorem, give more insight into the asymptotic relations stated
in Section 1.
\end{remark}

\begin{theorem} Large Gaussian quantiles satisfy the asymptotic representation
\[z_{1-\alpha}=\sqrt{-2\ln\alpha-\ln|\ln\alpha|-
\ln(4\pi)}\cdot
(1+O\frac{\ln|\ln\alpha|}{(\ln\alpha)^2}),
\alpha\rightarrow +0.\]
\end{theorem}

\subsection{Two-sided confidence intervals}
For $s\in\{1,2\}, \alpha>0$, put
$
L(s;\alpha)=\bar X_n-\frac{\sigma}{\sqrt{n}}z_{1-\alpha}(s)$
and
$ R(s;\alpha)=\bar X_n+\frac{\sigma}{\sqrt{n}}z^-_{1-\alpha}(s).$ Further, let
$ \alpha_i(n)>0,\,i=1,2, \alpha_1(n)+\alpha_2(n)<1,$ and
\[
ACI(s;\alpha_1(n),\alpha_2(n))=[L(s;\alpha_1(n)), R(s;\alpha_2(n))].
\]
If conditions $(\ref{Lgamma})$ and
$(\ref{Os1})$ are fulfilled
then
$P_\mu((-\infty, L(s;\alpha_1(n)))\; {\it covers}\; \mu)\sim \alpha_1(n)$ and $P_\mu(( R(s;\alpha_2(n)), \infty)\; {\it covers}\; \mu)\sim \alpha_2(n)$ as $n\rightarrow \infty$.

With more detailed notation
$\mu_{1,n}(s)=\mu_{1,n}(s;\alpha,\beta)$
and $\mu^-_{1,n}(s)=\mu^-_{1,n}(s;\alpha,\beta)$,

$P_{\mu_{1,n}(s;\alpha_1(n),\beta_1(n))}
((L(s;\alpha_1(n)),\infty)\; covers\; \mu)\sim \beta_1(n)$,

$P_{\mu^-_{1,n}(s;\alpha_2(n),\beta_2(n))}
((-\infty, R(s;\alpha_2(n)))\; covers\; \mu)\sim \beta_2(n), n\rightarrow \infty.$

The following corollary has thus been proved.
\begin{corollary}
 If $\alpha_1(n)\downarrow 0$, $\alpha_2(n)\downarrow 0$ as $n\rightarrow \infty$ and conditions ($\ref{Lgamma} $) and ($\ref{Os1}$)  are satisfied for $\gamma\in (\frac{1}{6},\frac{1}{4}]$ if $s=1$ and for  $\gamma\in (\frac{1}{4},\frac{3}{10}]$ if $s=2$, and with $(\alpha(n),\beta(n))=(\alpha_1(n),\alpha_2(n))$,
  then
 \[
 P_\mu (ACI(s;\alpha_1(n),\alpha_2(n)) {\; does\; not\; cover\; } \mu)\sim (\alpha_1(n)+\alpha_2(n)), \, n\rightarrow \infty.
 \]Moreover,
 \[\max\limits_{\nu\in\{\mu_{1,n}
 (s;\alpha_1(n),\beta_1(n)),
 \mu^-_{1,n}(s;\alpha_2(n),\beta_2(n)) \}}
 P_{\nu} (ACI(s) {\;  covers\; } \mu)\leq\max\{\beta_1(n),\beta_2(n)\}.
 \]
\end{corollary}

\section{Testing}
\subsection{Adjusted quantiles}
Let us consider the problem of testing the hypothesis
$H_0: \mu\leq\mu_0$ versus the alternative
$H_A: \mu > \mu_0$. The first and second kind adjusted  decision rules of the one-sided asymptotic Gauss test suggest to reject $H_0$ if
$T_{n,0}> z_{1-\alpha(n)}(s)$
 for $s=1$ or $s=2,$ respectively, where
 $T_{n,0}=\sqrt{n}(\bar X_n-\mu_0)/\sigma$. Because \[
  P_{\mu_0}(reject\; H_0)=P_{\mu_0}(ACI^u(s) \; does\; not\; cover\; \mu_0),\] it follows from Theorems 1 and 2 that under the conditions given there the (sequence of) probabilities of an error of first kind satisfy the asymptotic relation
  \[
  P_{\mu_0}(reject\; H_0)\sim \alpha(n), n\rightarrow \infty.
  \]Concerning the power function of this test, because
  \[
  P_{\mu_{1,n}(s)}(do\; not\; reject\; H_0)=  P_{\mu_{1,n}(s)}(ACI^u(s) \; covers\; \mu_0),
  \]it follows under the same assumptions that the probabilities of a second kind error in the case that the sequence of the modified local parameters is $(\mu_{1,n}(s))_{n=1,2,...}$,  satisfy
  \[
  P_{\mu_{1,n}(s)}(do\; not\; reject\; H_0)\sim \beta(n), n\rightarrow \infty.
  \]
Similar consequences for testing $H_0: \mu>\mu_0$ or $H_0: \mu\neq\mu_0$ are omitted, here.
\subsection{Adjusted statistics}
Let $T_n^{(1)}=T_n-\frac{g_1}{6\sqrt{n}}T_n^2$
and
$T_n^{(2)}=T_n^{(1)}-\frac{3g_2-8g_1^2}
{72n}T_n^3$ be the first and second kind adjusted asymptotically Gaussian statistics, respectively, where $T_n=\frac{\sqrt{n}}{\sigma}(\bar X_n - \mu)$.

\begin{theorem}
If the  conditions ($\ref{Lgamma} $) and ($\ref{Os1}$) are satisfied for a certain $\gamma \in (\frac{s}{2s+4},\frac{s+1}{2s+6}]$ where $s\in\{1,2\}$ then
\[
P_{\mu_0}(T_n^{(s)}>z_{1-\alpha(n)})\sim \alpha(n), \; n\rightarrow \infty
\]
 and
 \[
 P_{\mu_{1,n}(s)}(T_n^{(s)}\leq z_{1-\alpha(n)})\sim \beta (n), \;
n\rightarrow \infty. \]
\end{theorem}

Clearly, the results of this theorem apply to both hypothesis testing and confidence estimation in a similar way as described in the preceding sections.
\section{Application to exponential families} Let $\nu$ denote a $\sigma$-finite measure and assume that the distribution $P_\vartheta $ has the Radon-Nikodym density
$\frac{dP_\vartheta}{d\nu}(x)=
\frac{e^{\vartheta x}}{\int
e^{\vartheta x}\nu(dx)}=e^{\vartheta x-B(\vartheta)}$, say.
We assume that $X(\vartheta)\sim P_\vartheta$ and $X_1=X(\vartheta)-{ E}X(\vartheta)+\mu\sim\widetilde{P}_\mu$ where
$\vartheta$ is known and $\mu$ is unknown. In the product-shift-experiment
$ [R^n,{\cal B}^n,\;\{\widetilde{P}^{\times n}_\mu,\,\mu\in R\}]$, expectation estimation and testing may be done as in Sections 2 and 3, respectively, where
$g_1= B^{'''}(\vartheta)/(B^{''}(\vartheta))^{3/2}$
and $g_2$ allows a similar representation.

Another problem which can be dealt with is to test the hypothesis
$H_0:\vartheta\leq\vartheta_0$ versus the alternative $
H_{1n}:\vartheta\geq\vartheta_{1n}$  if one assumes that the expectation function
  $\vartheta\rightarrow B'(\vartheta)=E_\vartheta
X$ is strongly monotonous. For this case, we finally present just the following particular result which applies to both estimating and testing.

\begin{proposition}
If conditions ($\ref{Lgamma} $) and ($\ref{Os1}$) are satisfied for $\gamma\in(\frac{1}{6},\frac{1}{4}]$ then
\[P_{\vartheta_0}^{\times
n}(\sqrt{n}\frac{\overline{X}_n-B'(\vartheta_0)}{\sqrt{B''(\vartheta_0)}}>
z_{1-\alpha(n)}+\frac{B'''(\vartheta_0)}{6\sqrt{n}(B''(\vartheta_0))^{3/2}}
z^2_{1-\alpha(n)})\sim\alpha(n),n\;\rightarrow \infty.\]
\end{proposition}

\section{Sketch of proofs}
{\it Proof of Theorems 1 and 2.}
If condition ($\ref{Os1}$) is satisfied then $x=z_{1-\alpha(n)}=o(n^\gamma),n\rightarrow \infty$ for $\gamma \in (\frac{1}{6},\frac{3}{10}]$, and if ($\ref{Lgamma}$) then, according to \cite{Li},  $P_\mu(T_n>x)\sim f_{n,s}^{(X)}(x), x\rightarrow \infty$
where
$
f_{n,s}^{(X)}(x)=\frac{1}{\sqrt{2\pi}x}
\exp\{-\frac{x^2}{2}+\frac{x^3}
{\sqrt{n}}\sum\limits_{k=0}^{s-1}
a_k(\frac{x}{\sqrt{n}})^k\}
$ and $s$ is an integer satisfying
$\frac{s}{2(s+2)}<\gamma\leq\frac{s+1}{2(s+3)}$,
i.e. $s=1$ if $\gamma\in (\frac{1}{6}, \frac{1}{4}]$ and $s=2$ if $\gamma\in (\frac{1}{4}, \frac{3}{10}]$.
Here, $o(\cdot)$ stands for the small Landau symbol, and the constants $a_0=\frac{g_1}{6},\, a_1=\frac{g_2-3g_1^2}{24} $
are due to the skewness $g_1$ and kurtosis $g_2$ of $X$.
Note that $\frac{g_1x^2}{6\sqrt{n}}=o(x)$ because x$=o(n^{1/2})$, thus $x+\frac{g_1x^2}{6\sqrt{n}}=o(n^\gamma)$,
and $P_\mu(T_n>x+\frac{g_1x^2}{6\sqrt{n}})
\sim f_{n,1}(x+\frac{g_1x^2}{6\sqrt{n}})$.
Hence, $P_\mu(T_n>x+\frac{g_1x^2}{6\sqrt{n}})\sim 1-\Phi(x).$ Similarly,
$P_\mu(T_n>z_{1-\alpha(n)}(s))
\sim\alpha(n), \, s=1,2$. Further,
 $P_{\mu_{1,n}(s)}(T_n\leq z_{1-\alpha(n)}(s))$
\[=P_{\mu_{1,n}(s)}(\frac{\sqrt{n}}{\sigma}
 (\bar X_n -\mu_{1,n}(s))< z_{1-\alpha(n)}(s)
 -\frac{\sqrt{n}}{\sigma}(\mu_{1,n}(s)-\mu))
 =P_0(\frac{\sqrt{n}}{\sigma}
 \bar X_n< z_{\beta(n)}(s)).\]
The latter equality holds because $\{P_\mu, \mu\in (-\infty, \infty )\}$
  is assumed to be a shift family. It follows that
$P_{\mu_{1,n}(s)}(T_n\leq z_{1-\alpha(n)}(s))$
\[=P_0(\frac{\sqrt{n}}{\sigma}
 (-\bar X_n)\geq z_{1-\beta(n)}+
 \frac{-g_1}{6\sqrt{n}}z^2_{1-\beta(n)}
 +I_{\{2\}}(s)
 \frac{3g_2-4g_1}{72n}z^3_{1-\beta(n)}).\]

Note that $-g_1, g_2$ are skewness and kurtosis of $-X_1$. Thus,
\[P_{\mu_{1,n}(s)}(T_n\leq z_{1-\alpha(n)}(s))\sim f_{n,s}^{(-X)}(z_{1-\beta(n)}(s))
   \sim\beta(n), n \rightarrow \infty.\]
 Because
 $P_\mu(ACI^u \; does \; not \; cover\; \mu)=P_{\mu}(T_n > z_{1-\alpha(n)}(s))$ and

 $P_{\mu_{1,n}(s)}(ACI^u   \; covers\; \mu)
 =P_{\mu_{1,n}(s)}(T_n\leq z_{1-\alpha(n)}(s))$,
the theorems are proved.

\smallskip

{\it Proof of Remark 1.} The first statement of the remark follows from
\[
P_\mu(\mu>\bar X_n+\sigma z^{-}_{1-\alpha(n)}/\sqrt{n})
=P_\mu(\sqrt{n}(-\bar X_n+\mu)/\sigma >z^{-}_{1-\alpha(n)})
\]
and the second one from
\[
P_{\mu^{-}_{1,n}(s)}(\mu< \bar X_n+ \sigma z^{-}_{1-\alpha(n)}/\sqrt{n})
=P_0(\bar X_n>\mu-\mu^{-}_{1,n}(s)-\sigma z^{-}_{1-\alpha(n)}/\sqrt{n})\]
\[
=P_0(\sqrt{n}\bar X_n/\sigma > z_{1-\beta(n)}(s)).
\]

\smallskip

{\it Proof of Theorem 3.} We start from the well known relations
\[\alpha=1-\Phi(z_{1-\alpha})=
(1+O(\frac{1}{z^2_{1-\alpha}}))\frac{1}{\sqrt{2\pi}z_{1-\alpha}}
e^{-\frac{z^2_{1-\alpha}}{2}},\;\alpha\rightarrow 0.\]
The solution to the approximative
quantile equation
$\alpha=\frac{1}{\sqrt{2\pi}x}
e^{-\frac{x^2}{2}}$ will be denoted by $x=x_{1-\alpha}$.
Let us put
\begin{equation}
xe^{\frac{x^2}{2}}=\frac{1}{\sqrt{2\pi}\alpha}=:y.
\label{QGl}\end{equation}

If $x\geq 1$ then it follows from $(\ref{QGl})$ that $ y\geq e^{\frac{x^2}{2}}$, hence
$x^2\leq\ln(y^2)$. It follows again from $(\ref{QGl})$ that
$ y^2\leq\ln(y^2)e^{x^2}$, thus
$x^2\geq \ln(\frac{y^2}{\ln y^2}).$ After one more such step,
\[\ln(\frac{y^2}{\ln y^2})\leq
x^2\leq\ln[\frac{y^2}{\ln(\frac{y^2}{\ln y^2})}].\]
The theorem now follows from
\[x^2=\{\ln y^2-\ln 2-\ln\ln y\}\{1+O(\frac{\ln\ln y}{(\ln
y^2)^2})\}, y\rightarrow \infty.\]

\smallskip

{\it Proof of Theorem 4.} Recognize that if $g_{n,s}(x)=o(\frac{1}{x}), x\rightarrow \infty$ then $f^{(+/-)(X)}_{n,s}(x+g_{n,s}(x))\sim f^{(+/-)(X)}_{n,s}(x),x\rightarrow \infty.$ Let us restrict to the case
$s=1$. According to \cite{Li},
\[P_{\mu_0}(T_n^{(1)}>z_{1-\alpha(n)})\sim
P_{\mu_0}(\frac{3\sqrt{n}}{g_1}>T_n^{(1)}>z_{1-\alpha(n)}).\] The function $f^{(1)}_n(t)=t-\frac{g_1t^2}{6\sqrt{n}}$ has a positive derivative, $f^{(1)'}_n(t)=1-
\frac{g_1t}{3\sqrt{n}}>0$,
if $g_1t<3\sqrt{n}$. Denoting there the inverse function of $f_n^{(1)}$ by $f_n^{(1)^{-1}}$, it follows
$f_n^{(1)^{-1}}(x)=
x+\frac{g_1x^2}{6\sqrt{n}}+O(\frac{x^3}{n})$ and
$
f^{(1)}_n(f_n^{(1)^{-1}}(x)) =
x+o(\frac{1}{x}).
$ Thus,

$
P_{\mu_o}(T_n^{(1)}>z_{1-\alpha(n)})\sim
P_{\mu_o}(T_n>z_{1-\alpha(n)}(1))
\sim\alpha(n).$

Moreover,
$P_{\mu_{1n}(1;\alpha(n),\beta(n))}(T_n^{(1)}\leq
z_{1-\alpha(n)})=P_{\mu_{1n}(1;\alpha(n),\beta(n))}(T_n\leq T^{(1)^{-1}}_n
(z_{1-\alpha(n)}))$\\[1ex]
$=P_{\mu_{1n}(1;\alpha(n),\beta(n))}(\;\sqrt{n}\frac{\overline{X}_n
-\mu_{1n}(1)}{\sigma} \leq
z_{1-\alpha(n)}(1)+\frac{z^2_{1-\alpha(n)}g_1}{6\sqrt{n}}+O(\frac{z^3_{1-\alpha(n)}}{n})-
\sqrt{n}\frac{\mu_{1n}(1)-\mu_0}{\sigma}) \sim  f_{n,1}^{(-X)}(-z_{\beta(n)}(1)+ O(\frac{z^3_{1-\alpha(n)}}{n}))  \sim 1-\Phi(z_{1-\beta(n)})=\beta(n).$

\end{document}